\documentclass[12pt,a4article]{amsart}
\usepackage{amsmath,amsthm,amssymb,amscd}
\usepackage{algorithmic}

\newcommand{\Z}{\mathbb{Z}}
\newcommand{\T}{\mathcal{T}}

 \DeclareMathOperator{\id}{id}
 \DeclareMathOperator{\cyc}{c}
 \DeclareMathOperator{\dec}{d}
 \DeclareMathOperator{\len}{len}
 \DeclareMathOperator{\infs}{inf_s}
 \DeclareMathOperator{\sups}{sup_s}
 \DeclareMathOperator{\lens}{len_s}

\newtheorem{thm}{Theorem}[section]
\newtheorem{lem}[thm]{Lemma}
\newtheorem{prop}[thm]{Proposition}
\newtheorem{cor}[thm]{Corollary}
\newtheorem{alg}[thm]{Algorithm}

\theoremstyle{definition}

\newtheorem{test}{Test}
\theoremstyle{remark}
\newtheorem{rem}[thm]{Remark}

\begin{document}

\title{General cycling operations in Garside groups}
\author[H. Zheng]{Hao Zheng}
\address{Department of Mathematics, Zhongshan University, Guangzhou 510275, China}
\thanks{}
\email{zhenghao@mail.sysu.edu.cn}
\date{}

\begin{abstract}
In this article, we introduce the notion of cycling operations of
arbitrary order in Garside groups, which is a full generalization
of the cycling and decycling operations. Theoretically, this
notion together with other related concepts provides a context in
which various definitions and arguments concerning Garside groups
are unified and simplified as well as improved. Practically, it
yields a new algorithm which has a considerably improved
performance on solving the conjugacy problem of reducible braids.
\end{abstract}

\maketitle

\begin{small}

\end{small}

\section{Introduction}\label{sec:intro}

The solution of the conjugacy problem in braid groups backdated to
Garside \cite{Gar} who established the first algorithm to solve
the problem by means of calculating a conjugacy invariant of
braids, the so called summit set. In the past decade, with many
efforts (for example, \cite{EM, BKL1, BKL2, FG}) put on a refined
version of the summit set, the super summit set, the algorithm was
improved in various aspects. The algorithm and its improvements
were also applied to a large family of groups, known as the
Garside groups or small Gaussian groups \cite{DP}.

Recent progress on this issue was addressed to \cite{Geb}, in
which the super summit set was refined again to the ultra summit
set by posing the cycling-recurrence condition. Remarkably, the
algorithm resulted is so efficient that it makes practically
possible to solve the conjugacy problem of generic braids
(pseudo-Anosov braids) with large number of strands and word
length.

Nevertheless, in contrast with such success, when confined to a
specific class of braids (but still generic in practical sense),
the reducible braids, even the best algorithm due to \cite{Geb}
practically fails. We will justify this point in Section
\ref{sec:compute} by giving examples and experimental data. To sum
up, in the case of reducible braids, the cycling-recurrence
condition loses its control on the components, so the performance
of ultra summit set degenerates to the level of super summit set.

To remedy this deficiency, a natural way is to further refine the
ultra summit set by posing the cycling-recurrence condition on the
components of a reducible braid. At first sight, applying cycling
operation on the components requires knowledge of the reduction
system of a reducible braid. However, this is not the case. The
refinement is easily implemented by introducing the notion of
cycling operations of arbitrary order in Garside groups, which is
a full generalization of the cycling and decycling operations.
With a slight modification to the algorithms for computing super
summit set and ultra summit set, one is able to compute the fully
refined summit set effectively and achieve great performance
improvement on solving the conjugacy problem of reducible braids.

Apart from practical significance, the notion of the general
cycling operations turns out to be a very fundamental concept.
Together with the concepts of pushforward and pullback along
general cycling operations, it provides a context in which various
definitions and arguments concerning Garside groups are unified
and simplified, hence sheds light on these aspects (see Section
\ref{sec:application} and the end of Section \ref{sec:garside}).
From the theoretical point of view, these new concepts provide a
very convenient and powerful tool for future study of Garside
groups.

\tableofcontents

\section{Notations and basic facts}\label{sec:garside}

The notion of Garside group \cite{DP} is a natural generalization
of braid group and, more generally, Artin group of finite type. In
this section, we give a brief review of Garside groups and state
some basic facts and know results for later use or comparison.
Readers are referred to \cite{BKL1, DP, Pic, Deh, FG, Geb} for
more details.

Let $M$ be a monoid. We say $x \in M$ is an {\em atom} if $x \neq
1$ and $x = yz$ implies either $y=1$ or $z=1$. $M$ is said to be
{\em atomic} if it is generated by its atoms and for every $x \in
M$ there exists a finite number $\|x\|$, called the {\em norm} of
$x$, such that $x$ is a product of at most $\|x\|$ atoms.

A cancellative, atomic monoid $M$ is said to be {\em Gaussian} if
every two elements of $M$ have both a left (and right) greatest
common divisor and a left (and right) least common multiple.

A {\em Garside monoid} $M$ is a Gaussian monoid which admits a
Garside element. The {\em Garside element} is an element $\Delta
\in M$ such that its left divisors coincide with its right
divisors, they forming a finite set and generating $M$. The
divisors of the Garside element are called {\em simple elements}.

Every Garside monoid admits a group of fractions. A group $G$ is
called a {\em Garside group} if it is the group of fractions of a
Garside monoid.

The braid groups are main examples of Garside groups. If we write
the $n$-strand braid group in Artin presentation
$$B_n = \bigg\langle \sigma_1,\dots,\sigma_{n-1} \bigg|
  \begin{array}{ll}
    \sigma_i \sigma_j = \sigma_j \sigma_i, & |i-j| \geq 2 \\
    \sigma_i \sigma_j \sigma_i = \sigma_j \sigma_i \sigma_j, & |i-j| = 1
  \end{array} \bigg\rangle,
$$
then the monoid given by the same presentation is a Garside monoid
with Garside element (the half twist)
$$\Delta = (\sigma_1\cdots\sigma_{n-1}) (\sigma_1\cdots\sigma_{n-2}) \cdots (\sigma_1\sigma_2) \sigma_1.$$
This Garside structure of $B_n$ is referred to as the {\em
classical structure}.

An alternative Garside structure of $B_n$ was given in
\cite{BKL1}, referred to as the {\em dual structure} or {\em BKL
structure}. Since this structure will not be used in this article
we omit the precise description here.

Throughout this article, let $G$ denote a Garside group associated
with Garside monoid $M$ and Garside element $\Delta$. Let $S$ and
$A$ denote the finite sets of simple elements and atoms of $G$,
respectively.

For $x,y \in G$ we denote by $x \prec y$ the relation that $x$ is
a left divisor of $y$, i.e. $x^{-1}y \in M$, by $x \wedge y$ and
$x \vee y$ the left greatest common divisor and left least common
multiple of $x,y$ respectively. The conjugation $u^{-1}xu$ is
denoted as $x^u$ and the specific conjugation $\Delta^{-1} x
\Delta$ is also denoted as $\tau(x)$.

A fundamental fact about a Garside group $G$ is that, for every $x
\in G$, there is a unique decomposition $x = \Delta^p x_1 \cdots
x_l$, called the {\em (left) normal form} of $x$, satisfying the
conditions $x \wedge \Delta^{p+i} = \Delta ^p x_1 \cdots x_i$ and
$x_i \in S \setminus \{1,\Delta\}$. The {\em infimum}, {\em
supremum} and {\em canonical length} of $x$ are defined to be
$\inf x = p$, $\sup x = p+l$ and $\len x = l$, respectively.

The following basic facts will be repeatedly used in the article
without explanation.
\begin{enumerate}
\item The relation $\prec$ is a partial order.

\item $\tau(S)=S$. So $\tau^e = \id$ for some $e>0$ and $\Delta^e$
lies in the center of $G$.

\item $\tau(M)=M$. So $x \prec y$ if and only if $\tau(x) \prec
\tau(y)$. Moreover, we have $\tau(x \wedge y) = \tau(x) \wedge
\tau(y)$, $\tau(x \vee y) = \tau(x) \vee \tau(y)$ and $x\Delta
\wedge y\Delta = (x \wedge y)\Delta$, $x\Delta \vee y\Delta = (x
\vee y)\Delta$.

\item $\Delta^{\inf x} \prec x \prec \Delta^{\sup x}$. So, $1
\prec x^n\Delta^{-n\inf x}$ and $1 \prec x^{-n}\Delta^{n\sup x}$
for $n\geq0$.

\item Each set $\{ x \mid p_1 \leq \inf x, \; \sup x \leq p_2 \}$
is finite. So, a sequence of sufficient length in it has
repetitions.
\end{enumerate}

The approach found by Garside to solve the conjugacy problem in a
Garside group $G$ is to associate to each element $x \in G$ a
computable, nonempty subset $\tilde{C}(x) \subset C(x)$ which is
only dependent on the conjugacy class $C(x)$ of $x$. Then given
two elements $x,y \in G$, one just computes and compares
$\tilde{C}(x)$ and $\tilde{C}(y)$ to see whether $x$ and $y$
belong to the same conjugacy class. The summit set and its
refinements, the super summit set and the ultra summit set, are
such type of conjugacy invariants.

The {\em summit infimum}, {\em summit supremum} and {\em summit
length} of the conjugacy class $C(x)$ of $x$ are
\begin{align*}
  & \infs x = \max \{ \inf y \mid y \in C(x) \}, \\
  & \sups x = \min \{ \sup y \mid y \in C(x) \}, \\
  & \lens x = \sups x - \infs x.
\end{align*}
The {\em cycling} and {\em decycling} operations on $x = \Delta^p
x_1 \cdots x_l$ in normal form are the conjugations
\begin{align*}
  & \cyc(x) = x^{\Delta^px_1\Delta^{-p}} = \Delta^p x_2 \cdots x_l \tau^{-p}(x_1), \\
  & \dec(x) = x^{\Delta^p x_1 \cdots x_{l-1}} = \Delta^p \tau^p(x_l) x_1 \cdots x_{l-1}.
\end{align*}
Note that both operations neither decrease the infimum nor
increase the supremum.

With these notations, one defines the {\em super summit set}
$$C^s(x) = \{ y \in C(x) \mid \inf y = \infs x, \;\; \sup y = \sups x \},$$
and the {\em ultra summit set}
$$C^u(x) = \{ y \in C^s(x) \mid \text{$\cyc^N(y) = y$ for some $N>0$} \}.$$
The finiteness of both sets are clear. The nonemptiness and the
computability of these conjugacy invariants can be derived from
the following theorems (see the references linked).

\begin{thm}[\cite{EM,BKL1,Pic}]\label{thm:a0}
If $\inf x < \infs x$ then $\inf\cyc^N(x) > \inf x$ for some
$N>0$. Similarly, if $\sup x > \sups x$ then $\sup\dec^N(x) < \sup
x$ for some $N>0$.
\end{thm}

\begin{thm}[\cite{FG}]\label{thm:a1}
If $x^u, x^v \in C^s(x)$ then $x^{u \wedge v} \in C^s(x)$.
\end{thm}

\begin{thm}[\cite{Geb}]\label{thm:a2}
If $x^u, x^v \in C^u(x)$ then $x^{u \wedge v} \in C^u(x)$.
\end{thm}

As an evidence of the powerfulness of our new machinery, all these
theorems will appear as easy corollaries in the next section.

\section{New definitions and main results}\label{sec:cycling}

The {\em cycling operation of order $q$} on $x$ is the conjugation
$$\cyc_q(x) = x^{x \wedge \Delta^q}.$$
We have the $\cyc_q$-recurrence set
$$G_q = \{ x \in G \mid \text{$\cyc_q^N(x) = x$ for some $N>0$} \}.$$
The following properties are immediate  from definition. In
particular, the last one says that each sequence $x, \cyc_q(x),
\cyc_q^2(x), \dots$ eventually runs into a closed orbit, so $C(x)
\cap G_q$ is always nonempty.

\begin{lem}\label{lem:A1}
Properties of cycling operations.
\begin{enumerate}
  \item $\tau \cyc_q(x) = \cyc_q \tau(x)$. So, $\tau(G_q) = G_q$.
  \item For $x = \Delta^p x_1 \cdots x_l$ in normal form, we have
  $$\cyc_q(x) = \left\{ \begin{array}{ll}
       \tau^q(x), & q \leq \inf x, \\
       x_{q-p+1} \cdots x_l \Delta^p x_1 \cdots x_{q-p}, & \inf x < q < \sup x, \\
       x, & q \geq \sup x. \\
    \end{array} \right.
  $$
  \item $x \in G_q$ for $q\leq\inf x$ or $q\geq\sup x$.
  \item $\inf x \leq \inf \cyc_q(x)$ and $\sup \cyc_q(x) \leq \sup x$.
\end{enumerate}
\end{lem}

These new cycling operations are indeed natural generalizations of
the cycling and decycling operations. Note that
$$\cyc(x) = \tau^{-\inf x} \cyc_{\inf x+1}(x), \quad
  \dec(x) = \cyc_{\sup x-1}(x).
$$
In the next section we derive the following theorem. Since for
each $x$ the $\cyc_q$ orbit of $x$ eventually runs into $G_q$, it
follows from the theorem that if $\inf x < q \leq \infs x$ then
$\inf\cyc_q^N(x) \geq q$ holds for sufficient large $N$ and a
similar statement for supremum. In particular, the specific case
$q=\inf x+1$ or $q=\sup x-1$ gives rise to Theorem \ref{thm:a0}.

\begin{thm}\label{thm:A0}
We have $x \not\in G_q$ for $\inf x < q \leq \infs x$ or $\sups x
\leq q < \sup x$.
\end{thm}

In the sequel, the super summit set and the ultra summit set are
nothing but
\begin{align*}
  C^s(x) & = C(x) \cap \bigcap_{q \in \{\infs x,\sups x\}} G_q, \\
  C^u(x) & = C(x) \cap \bigcap_{q \in \{\infs x,\infs x+1,\sups x\}} G_q.
\end{align*}
The fully refined summit set we define here is
$$C^*(x) = C(x) \cap \bigcap_{q \in \Z} G_q
  = C(x) \cap \bigcap_{\infs x \leq q \leq \sups x} G_q.
$$
Remark the obvious inclusions
$$C^*(x) \subset C^u(x) \subset C^s(x).$$

\smallskip

The following theorem will also be proved in the next section. As
an immediate consequence, we conclude that $x^{u \wedge v} \in
C^*(x)$ (resp. $C^s(x), C^u(x)$) provided $x^u, x^v \in C^*(x)$
(resp. $C^s(x), C^u(x)$). So, we reach an alternative proof of
Theorem \ref{thm:a1} and Theorem \ref{thm:a2}.

\begin{thm}\label{thm:A1}
If $x^u, x^v \in G_q$ then $x^{u \wedge v} \in G_q$. In
particular, $\cyc_q(G_p) \subset G_p$ for all $p,q \in \Z$.
\end{thm}

From Lemma \ref{lem:A1}(3) and the inclusion $\cyc_q(G_p) \subset
G_p$, we have the following algorithm. In particular, the set
$C^*(x)$ is always nonempty.

\begin{alg}\label{alg:A1}
Given an element $x$ of $G$, the following algorithm computes an
element of $C^*(x)$.
\begin{algorithmic}\upshape
  \STATE Set $q=\inf x+1$.
  \WHILE {$q < \sup x$}
    \STATE Compute $x, \cyc_q(x), \cyc_q^2(x), \dots, \cyc_q^N(x)$
    until repetition encountered.
    \STATE Set $x=\cyc_q^N(x)$.
    \STATE Set $q=q+1$.
  \ENDWHILE
  \STATE \textbf{return} $x$
\end{algorithmic}
\end{alg}

Now we proceed to present an algorithm for computing the whole
$C^*(x)$. Define the {\em full cycling trajectory} of $x \in G$
$$T(x) = \{ \cyc_{q_k}\cdots\cyc_{q_2}\cyc_{q_1}(x) \mid q_i \in \Z, \; k \geq 0 \}.$$
The validity of the next algorithm follows from Lemma
\ref{lem:A1}(1),(2).

\begin{alg}\label{alg:A2}
Given an element $x$ of $G$, the following algorithm computes the
full cycling trajectory $T(x)$.
\begin{algorithmic}\upshape
  \STATE Set $T = \{ x, \tau(x), \dots, \tau^{e-1}(x) \}$
  where $e>0$ satisfies $\tau^e=\id$.
  \FOR {$y \in T$}
    \STATE Set $T = T \cup \{ \cyc_q(y) \mid \inf y < q < \sup y \}$.
  \ENDFOR
  \STATE \textbf{return} $T$
\end{algorithmic}
\end{alg}

Let $A(x)$ denote the set of $\prec$-minimal elements in $\{ u \in
S \setminus \{1\} \mid x^u \in C^*(x) \}$. The following theorems
are proved in the next section. Thanks to them we have Algorithm
\ref{alg:A3} for computing $C^*(x)$.

\begin{thm}\label{thm:A2}
For each pair $x_1, x_2 \in C^*(x)$ there exists a sequence
$$y_1=x_1, \; y_2, \; y_3, \dots, \; y_k=x_2 \in C^*(x)$$
such that $y_{i+1}=y_i^{u_i}$ for some $u_i \in A(y_i)$.
\end{thm}

\begin{thm}\label{thm:A3}
For each pair $x_1, x_2 \in C^*(x)$ with $T(x_1)=T(x_2)$ and for
each $u_1 \in A(x_1)$, there exists $u_2 \in A(x_2)$ such that
$T(x_1^{u_1}) = T(x_2^{u_2})$.
\end{thm}

\begin{alg}\label{alg:A3}
Given an element $x$ of $G$, the following algorithm computes
$C^*(x)$.
\begin{algorithmic}\upshape
  \STATE Compute $\tilde{x} \in C^*(x)$ and set $\T = \{ T(\tilde{x}) \}$.
  \FOR {$T \in \T$}
    \STATE Choose $y \in T$ and set $\T = \T \cup \{ T(y^u) \mid u \in A(y) \}$.
  \ENDFOR
  \STATE \textbf{return} $\bigcup_{T \in \T} T$
\end{algorithmic}
\end{alg}

Note that Algorithm \ref{alg:A3} involves a computation of the set
$\{ T(y^u) \mid u \in A(y) \}$, which we will work out in Section
\ref{sec:push_pull2} along the lines of \cite{FG, Geb}. Although
we can alternatively compute the superset $\{ T(y^u) \mid u \in S
\setminus \{1\}, \; y^u \in C^*(x) \}$ in the algorithm, which is
much easier to be implemented, as argued in \cite{FG} this may
decrease the performance considerably, because a Garside group may
have a large number of simple elements while only a few atoms. For
example, the braid group $B_n$ endowed with the classical Garside
structure has $n!$ simple elements but only $n-1$ atoms. So a
delicate implementation of Algorithm \ref{alg:A3} is necessary for
practical use.

\section{Pushforward and pullback I}\label{sec:push_pull1}

The notions of pushforward and pullback were introduced in
\cite{Geb} (pushforward was called transport instead) where they
were used to keep track of the cycling orbits of various
conjugations of an element $x$ and were proved to be very powerful
in the study of ultra summit set.

These notions are also applicable for general cycling operations.
Inspiringly, in this new setting they can be defined in a very
concise form. The {\em pushforward} $\phi_{x,q}(u)$ and {\em
pullback} $\pi_{x,q}(u)$ of $u$ along the cycling operation $x \to
\cyc_q(x)$ are defined as
\begin{align*}
  & \phi_{x,q}(u) = x''u \wedge x'^{-1}\Delta^q\tau^q(u), \\
  & \pi_{x,q}(u) = \Delta^{\inf u} \vee x''^{-1}u \vee x'\Delta^{-q}\tau^{-q}(u),
\end{align*}
respectively, where $x'=x\wedge\Delta^q$ and $x''=x'^{-1}x$. We
clarify these definitions by a pair of lemmas.

In what follows, if the context is clear we omit the subscripts of
$\phi, \pi$.

\begin{lem}\label{lem:B1}
Properties of pushforward.
\begin{enumerate}
\item $(x\wedge\Delta^q)\phi(u) = u (x^u\wedge\Delta^q)$, so
$\cyc_q(x)^{\phi(u)} = \cyc_q(x^u)$. See the diagram below.
$$\begin{CD}
    x^u @>{x^u\wedge\Delta^q}>> \cyc_q(x^u) \\
    @A{u}AA @AA{\phi(u)}A \\
    x @>{x\wedge\Delta^q}>> \cyc_q(x)
  \end{CD}
$$
\item $\phi(\Delta^p) = \Delta^p$.

\item If $u \prec v$ then $\phi(u) \prec \phi(v)$.

\item $\inf u \leq \inf\phi(u)$ and $\sup\phi(u) \leq \inf u$.

\item $\phi(u \wedge v) = \phi(u) \wedge \phi(v)$.

\item If $x^u=x^v$ and $\phi(u)=\phi(v)$ then $u=v$.
\end{enumerate}
\end{lem}

\begin{proof}
(1) $(x\wedge\Delta^q)\phi(u) = x'(x''u \wedge x'^{-1}u\Delta^q) =
xu \wedge u\Delta^q = u (x^u\wedge\Delta^q)$.

(2) $\phi(\Delta^p) = x''\Delta^p \wedge x'^{-1}\Delta^q\Delta^p =
x'^{-1} (x \wedge \Delta^q) \Delta^p = \Delta^p$.

(3) If $u \prec v$ then $\phi(u) = x''u \wedge
x'^{-1}\Delta^q\tau^q(u) \prec x''v \wedge
x'^{-1}\Delta^q\tau^q(v) = \phi(v)$.

(4) Since $\Delta^{\inf u} \prec u \prec \Delta^{\sup u}$, from
(2) and (3) we have $\Delta^{\inf u} = \phi(\Delta^{\inf u}) \prec
\phi(u) \prec \phi(\Delta^{\sup u}) = \Delta^{\sup u}$. Hence
$\inf u \leq \inf\phi(u)$ and $\sup\phi(u) \leq \sup u$.

(5) $\phi(u \wedge v) = x''(u \wedge v) \wedge
x'^{-1}\Delta^q\tau^q(u \wedge v) = x''u \wedge x''v \wedge
x'^{-1}\Delta^q\tau^q(u) \wedge x'^{-1}\Delta^q\tau^q(v) = \phi(u)
\wedge \phi(v)$.

(6) By (1) we have $u = (x\wedge\Delta^q) \phi(u)
(x^u\wedge\Delta^q)^{-1}$ and $v = (x\wedge\Delta^q) \phi(v)
(x^v\wedge\Delta^q)^{-1}$. So $x^u=x^v$ and $\phi(u)=\phi(v)$
imply $u=v$.
\end{proof}

\begin{lem}\label{lem:B2}
Properties of pullback.
\begin{enumerate}
\item $\pi(u) \prec v$ if and only if $\inf u \leq \inf v$ and $u
\prec \phi(v)$.

\item $\pi(\Delta^p) = \Delta^p$.

\item If $u \prec v$ then $\pi(u) \prec \pi(v)$.

\item $\inf u \leq \inf\pi(u)$ and $\sup\pi(u) \leq \inf u$.

\item $u \prec \phi\pi(u)$.

\item If $\inf \phi(v)= \inf v$ then $\pi\phi(v) \prec v$.
\end{enumerate}
\end{lem}

\begin{proof}
(1) By the definition of pushforward, we have $u \prec \phi(v)
\Longleftrightarrow u \prec x''v$ and $u \prec
x'^{-1}\Delta^q\tau^q(v) \Longleftrightarrow x''^{-1}u \prec v$
and $x'\Delta^{-q}\tau^{-q}(u) \prec v \Longleftrightarrow
x''^{-1}u \vee x'\Delta^{-q}\tau^{-q}(u) \prec v$. Therefore,
$\pi(u) \prec v$ if and only if $\Delta^{\inf u} \prec v$ and $u
\prec \phi(v)$.

(2) $\pi(\Delta^p) = \Delta^p \vee x''^{-1}\Delta^p \vee
x'\Delta^{-q}\Delta^p = (1 \vee x''^{-1} \vee x'\Delta^{-q})
\Delta^p = \Delta^p$.

(3) If $u \prec v$ then $\pi(u) = \Delta^{\inf u} \vee x''^{-1}u
\vee x'\Delta^{-q}\tau^{-q}(u) \prec \Delta^{\inf v} \vee
x''^{-1}v \vee x'\Delta^{-q}\tau^{-q}(v) = \pi(v)$.

(4) Follows from (2) and (3).

(5),(6) Apply (1) for $v = \pi(u)$ and $u = \phi(v)$ respectively.
\end{proof}

For specific values of $q$, we have several more properties of the
pushforward.

\begin{lem}\label{lem:B3}
For $q \leq \inf x$ we have
\begin{enumerate}
\item $\tau^{-q}\cyc_q(x) = x$ and $\tau^{-q}\phi(u) =
u(x^u\wedge\Delta^q)\Delta^{-q}$, so $x^{\tau^{-q}\phi(u)} =
\tau^{-q}\cyc_q(x^u)$;
$$\begin{CD}
    x^u @>{x^u\wedge\Delta^q}>> \cyc_q(x^u) @>{\Delta^{-q}}>> \tau^{-q}\cyc_q(x^u) \\
    @A{u}AA @AA{\phi(u)}A @AA{\tau^{-q}\phi(u)}A \\
    x @>{\Delta^q}>> \cyc_q(x) @>{\Delta^{-q}}>> x \\
  \end{CD}
$$

\item $\tau^{-q}\phi(u) \prec u$ with equality holds if and only
if $\inf x^u \geq q$.
\end{enumerate}
Similarly, For $q \geq \sup x$ we have
\begin{enumerate}
\item $\cyc_q(x) = x$ and $\phi(u) = u(x^u\wedge\Delta^q)$, so
$x^{\phi(u)} = \cyc_q(x^u)$;

\item $\phi(u) \prec u$ with equality holds if and only if $\sup
x^u \leq q$.
\end{enumerate}
\end{lem}

\begin{proof}
For $q \leq \inf x$ we have $x\wedge\Delta^q=\Delta^q$. Then (1)
is a special case of Lemma \ref{lem:B1}(1). Moreover, from
$\tau^{-q}\phi(u) = u(x^u\wedge\Delta^q)\Delta^{-q}$ we have
$\tau^{-q}\phi(u) \prec u\Delta^q\Delta^{-q} = u$ with equality
holds if and only if $\Delta^q \prec x^u$. Hence (2) holds.

The supremum part is proved similarly.
\end{proof}

Combining pushforwards and pullbacks along single cycling
operations, one has the pushforward and pullback along an
arbitrary cycling orbit
$$x \to \cyc_{q_1}(x) \to \cyc_{q_2}\cyc_{q_1}(x) \to \cdots \to \cyc_{q_k}\cdots\cyc_{q_1}(x).$$
For example, the pushforward $\phi_{x,q}^{(n)}(u)$ and pullback
$\pi_{x,q}^{(n)}(u)$ of $u$ along the cycling orbit $x \to
\cyc_q(x) \to \cyc_q^2(x) \to \cdots \cyc_q^n(x)$ are defined by
induction as
\begin{align*}
  & \phi_{x,q}^{(0)}(u) = u, \quad
  \phi_{x,q}^{(n)}(u) = \phi_{\cyc_q^{n-1}(x),q}\phi_{x,q}^{(n-1)}(u), \\
  & \pi_{x,q}^{(0)}(u) = u, \quad
  \pi_{x,q}^{(n)}(u) = \pi_{x,q}\pi_{\cyc_q(x),q}^{(n-1)}(u).
\end{align*}
Now suppose $x \in G_q$ and let $L$ be the $\cyc_q$-orbit length
of $x$, i.e. the minimal positive integer such that
$\cyc_q^L(x)=x$. The pushforward and pullback of $u$ along the
cycling orbit $x \to \cyc_q(x) \to \cyc_q^2(x) \to \cdots \to
\cyc_q^L(x)$ will be denoted as
\begin{align*}
  & \tilde\phi_{x,q}(u) = \phi_{x,q}^{(L)}(u), \quad \quad
  \tilde\pi_{x,q}(u) = \pi_{x,q}^{(L)}(u).
\end{align*}

\medskip

The following proposition plays a crucial role in this article.
Thanks to it, all theorems we claimed in the previous section are
derived readily.

\begin{prop}\label{prop:B1}
Suppose $x \in G_q$. Then $x^u \in G_q$ if and only if
$\tilde\phi^N(u) = u$ for some $N>0$.
\end{prop}

\begin{proof}
Let $L$ be the $\cyc_q$-orbit length of $x$. If $\tilde\phi^N(u) =
u$ then by Lemma \ref{lem:B1}(1) $\cyc_q^{NL}(x^u) =
\cyc_q^{NL}(x)^{\tilde\phi^N(u)} = x^u$ thus $x^u \in G_q$.

Conversely, suppose $x^u \in G_q$. By Lemma \ref{lem:B1}(4) the
equality $\tilde\phi^{N_1}(u) = \tilde\phi^{N_2}(u)$ holds for
some $0 \leq N_1<N_2$. Then from Lemma \ref{lem:B1}(1) we have
$$\cyc_q^{N_1L}(x^u) = x^{\phi^{(N_1L)}(u)} = x^{\phi^{(N_2L)}(u)} = \cyc_q^{N_2L}(x^u).$$
Therefore, $\cyc_q^{N_1L-i}(x^u) = \cyc_q^{N_2L-i}(x^u)$ for
$i=0,1,\dots,N_1L$, since $x^u \in G_q$. By Lemma \ref{lem:B1}(1)
again
\begin{align*}
  \cyc_q^{N_1L-i}(x)^{\phi^{(N_1L-i)}(u)}
  & = \cyc_q^{N_1L-i}(x^u) \\
  & = \cyc_q^{N_2L-i}(x^u)
  = \cyc_q^{N_2L-i}(x)^{\phi^{(N_2L-i)}(u)}
  = \cyc_q^{N_1L-i}(x)^{\phi^{(N_2L-i)}(u)}
\end{align*}
for $i=0,1,\dots,N_1L$. Note that in the last equality we used
$\cyc_q^L(x)=x$. Finally, starting from $\phi^{(N_1L)}(u) =
\phi^{(N_2L)}(u)$ and applying Lemma \ref{lem:B1}(6) $N_1L$ times
we get $u=\phi^{(N_2L-N_1L)}(u)$, i.e.
$u=\tilde\phi^{N_2-N_1}(u)$.
\end{proof}

\begin{proof}[Proof of Theorem \ref{thm:A0}]
Suppose $x \in G_q$ for some $\inf x < q \leq \infs x$. Choose $y,
u$ such that $\inf y \geq q$ and $x = y^u$. By Lemma
\ref{lem:B3}(2) we have
$$(\tau^{-q}\phi_{y,q})^N(u)
  \prec (\tau^{-q}\phi_{y,q})^{N-1}(u)
  \prec \cdots \prec \tau^{-q}\phi_{y,q}(u) \prec u.
$$
Since both $y,y^u \in G_q$, using the same argument in the proof
of Proposition \ref{prop:B1} we can show
$(\tau^{-q}\phi_{y,q})^N(u) = u$ for some $N>0$. It follows that
$\tau^{-q}\phi_{y,q}(u) = u$. By Lemma \ref{lem:B3}(2) we have
$\inf y^u \geq q$, which contradicts the assumption $\inf y^u =
\inf x < q$.

The other case of the theorem is proved similarly.
\end{proof}

\begin{proof}[Proof of Theorem \ref{thm:A1}]
For $x \in G_q$, if $x^u, x^v \in G_q$ then by Proposition
\ref{prop:B1} there exists $N>0$ such that $\tilde\phi^N(u)=u$ and
$\tilde\phi^N(v)=v$. So by Lemma \ref{lem:B1}(5) $\tilde\phi^N(u
\wedge v) = \tilde\phi^N(u) \wedge \tilde\phi^N(v) = u \wedge v$.
Applying Proposition \ref{prop:B1} again yields $x^{u \wedge v}
\in G_q$.

For general $x \in G$, suppose $x^w \in G_q$. Then
$(x^{w})^{w^{-1}u}, (x^{w})^{w^{-1}v} \in G_q$. Therefore $x^{u
\wedge v} = (x^{w})^{w^{-1}u \wedge w^{-1}v} \in G_q$.

As to the inclusion $\cyc_q(G_p) \subset G_p$, one notices that if
$x \in G_p$ then $x^x, x^{\Delta^q} \in G_p$ hence $\cyc_q(x) =
x^{x \wedge \Delta^q} \in G_p$.
\end{proof}

\begin{cor}\label{cor:B1}
For $x \in G_q$, $\phi_{x,q}$ restricts to a bijection
$$\phi_{x,q} : \{ u \mid x^u \in C^*(x) \} \to \{ u \mid \cyc_q(x)^u \in C^*(x) \}.$$
\end{cor}

\begin{proof}
For $x^u \in C^*(x)$, it follows from Lemma \ref{lem:B1}(1) and
the inclusion $\cyc_q(G_p) \subset G_p$ that
$\cyc_q(x)^{\phi_{x,q}(u)} = \cyc_q(x^u) \in C^*(x)$. Therefore,
$\phi_{x,q}$ does restrict to above map. By Proposition
\ref{prop:B1}, the map is invertible.
\end{proof}

\begin{proof}[Proof of Theorem \ref{thm:A2}]
Suppose $x_2 = x_1^u$. Multiplying a power of $\Delta$ if
necessary, we assume $u \in M$. The theorem is proved by induction
on the norm $\|u\|$. First, set $y_1=x_1$. If $\|u\|=0$ then
$x_1=x_2$ and we have nothing to do. Otherwise, since $x_1^u,
x_1^\Delta \in C^*(x)$, by Theorem \ref{thm:A1} $x_1^{u \wedge
\Delta} \in C^*(x)$ hence we can choose $u_1 \in A(x_1)$ such that
$u_1 \prec u \wedge \Delta$. Set $y_2=y_1^{u_1}$. Then
$y_2^{u_1^{-1}u} = x_2$ but $\|u_1^{-1}u\| < \|u\|$. By inductive
hypothesis, there exists a sequence $y_2, \; y_3, \; \dots, \;
y_k=x_2$ such that $y_{i+1}=y_i^{u_i}$ for some $u_i \in A(y_i)$.
\end{proof}

\begin{proof}[Proof of Theorem \ref{thm:A3}]
Suppose $x_1, x_2$ are connected by a cycling orbit
$$x_1 \to \cyc_{q_1}(x_1) \to \cyc_{q_2}\cyc_{q_1}(x_1) \to \cdots \to \cyc_{q_k}\cdots\cyc_{q_2}\cyc_{q_1}(x_1)=x_2.$$
By Corollary \ref{cor:B1}, the pushforward $\psi$ along the orbit
restricts to a bijection
$$\psi : \{ u \mid x_1^u \in C^*(x) \} \to \{ u \mid x_2^u \in C^*(x) \}$$
By Lemma \ref{lem:B1}(2),(3), $\psi$ further restricts to a
bijection $\psi : A(x_1) \to A(x_2)$. From Lemma \ref{lem:B1}(1),
we have $x_2^{\psi(u_1)} =
\cyc_{q_k}\cdots\cyc_{q_2}\cyc_{q_1}(x_1^{u_1})$ thus
$T(x_2^{\psi(u_1)}) = T(x_1^{u_1})$ for each $u_1 \in A(x_1)$.
\end{proof}

Till now, we have not made any use of the notion of pullback. We
conclude this section by preparing the following proposition.
Roughly speaking, $\tilde\pi$-recurrency guarantees lower bounds
for $\tilde\phi$-orbits (see also Lemma \ref{lem:C2}(1)).

\begin{prop}\label{prop:B2}
Suppose $x \in G_q$ and $\tilde\pi^{N_1}(u) = \tilde\pi^{N_2}(u)$
for some $0 \leq N_1 < N_2$. Then for every $v$ satisfying
$\tilde\pi^{N_1}(u) \prec v$ there exists arbitrarily large $N$
such that $u \prec \tilde\phi^N(v)$.
\end{prop}

\begin{proof}
From Lemma \ref{lem:B2}(5) and the hypotheses $\tilde\pi^{N_1}(u)
= \tilde\pi^{N_2}(u)$, $\tilde\pi^{N_1}(u) \prec v$, we have
$$u \prec \tilde\phi^{N_1+k(N_2-N_1)}\tilde\pi^{N_1+k(N_2-N_1)}(u)
  = \tilde\phi^{N_1+k(N_2-N_1)}\tilde\pi^{N_1}(u)
  \prec \tilde\phi^{N_1+k(N_2-N_1)}(v)$$
for any $k\geq0$.
\end{proof}

\section{Pushforward and pullback II}\label{sec:push_pull2}

In this section we derive a delicate implementation of Algorithm
\ref{alg:A3}.

For $x \in \cap_{q\in\Z}G_q$, let $\Phi_x$ be the free group
generated by $\{ \tilde\phi_{x,q} \mid q \in \Z \}$ which, by
Corollary \ref{cor:B1}, acts on the set $\{ u \mid x^u \in C^*(x)
\}$.

\begin{lem}\label{lem:C1}
Properties of $\Phi_x$.
\begin{enumerate}
\item $T(x^{\psi(u)}) = T(x^u)$ for every $\psi\in\Phi_x$.

\item $\Phi_x$ preserves the partial order $\prec$.

\item $\inf u = \inf \psi(u)$ and $\sup \psi(u) = \sup u$ for
every $\psi\in\Phi_x$.

\item If $u \prec v$ and $\Phi_xu = \Phi_xv$ then $u=v$.
\end{enumerate}
\end{lem}

\begin{proof}
(1),(2),(3) Follows from Lemma \ref{lem:B1}(1),(3),(4)
respectively, together with Proposition \ref{prop:B1}.

(4) Suppose $v=\psi(u)$ for some $\psi\in\Phi_x$. From (2) we have
$u \prec \psi(u) \prec \psi^2(u) \prec \cdots$. But by (3)
$\psi^N(u) = u$ holds for some $N>0$. It follows that $u =
\psi(u)$, hence $u = v$.
\end{proof}

\begin{prop}\label{prop:C1}
The relation $\prec$ defined on the set of $\Phi_x$-orbits by
$$\text{$\Phi_xu \prec \Phi_xv$ if $\psi_1(u) \prec \psi_2(v)$ for some $\psi_1,\psi_2\in\Phi_x$}$$
is a partial order.
\end{prop}

\begin{proof}
Reflexivity and transitivity are clear. Suppose $\Phi_xu \prec
\Phi_xv \prec \Phi_xu$, i.e. $\psi_1(u) \prec \psi_2(v)$ and
$\psi_1'(v) \prec \psi_2'(u)$ for some
$\psi_1,\psi_2,\psi_1',\psi_2'\in\Phi_x$. Then by Lemma
\ref{lem:C1}(2), $u \prec \psi_1^{-1}\psi_2(v) \prec
\psi_1^{-1}\psi_2\psi_1'^{-1}\psi_2'(u)$ and by Lemma
\ref{lem:C1}(4), $u = \psi_1^{-1}\psi_2\psi_1'^{-1}\psi_2'(u)$.
Hence $u=\psi_1^{-1}\psi_2(v)$ and $\Phi_xu = \Phi_xv$, i.e. the
relation is symmetric.
\end{proof}

Thanks to Theorem \ref{thm:A1} we have a well defined map
$$\mu_x : G \to \{ v \mid x^v \in C^*(x) \}, \quad
  u \mapsto \wedge \{ v \mid u \prec v, \; x^v \in C^*(x) \}.
$$
That is, $\mu_x(u)$ is the $\prec$-minimal element satisfying $u
\prec \mu_x(u)$ and $x^{\mu_x(u)} \in C^*(x)$. Slightly abusing
notation, we denote by $\Phi_xu$ the $\Phi_x$-orbit of $\mu_x(u)$
for arbitrary $u \in G$.

By definition, $\mu_x(u) \prec \mu_x(v)$ and $\Phi_xu \prec
\Phi_xv$ hold for $u \prec v$. This fact will be used repeatedly
in the remainder of this section.

\begin{lem}\label{lem:C2}
For $x \in \cap_{q\in\Z}G_q$ we have the followings.
\begin{enumerate}
\item $\Phi_x\tilde\phi_{x,q}(u) \prec \Phi_xu$.

\item $\Phi_x\tilde\pi_{x,q}(u) \prec \Phi_xu$.
\end{enumerate}
\end{lem}

\begin{proof}
(1) $\Phi_x\tilde\phi_{x,q}(u) \prec
\Phi_x\tilde\phi_{x,q}\mu_x(u) = \Phi_xu$.

(2) By Lemma \ref{lem:C1}(3), $\inf
\tilde\phi_{x,q}\tilde\phi_{x,q}^{-1}\mu_x(u) = \inf
\tilde\phi_{x,q}^{-1}\mu_x(u)$. Applying Lemma \ref{lem:B2}(6)
yields
$\tilde\pi_{x,q}\tilde\phi_{x,q}\tilde\phi_{x,q}^{-1}\mu_x(u)
\prec \tilde\phi_{x,q}^{-1}\mu_x(u)$. Therefore,
$$\Phi_x\tilde\pi_{x,q}(u)
  \prec \Phi_x\tilde\pi_{x,q}\mu_x(u)
  = \Phi_x\tilde\pi_{x,q}\tilde\phi_{x,q}\tilde\phi_{x,q}^{-1}\mu_x(u)
  \prec \Phi_x\tilde\phi_{x,q}^{-1}\mu_x(u)
  = \Phi_xu.
$$
\end{proof}

\begin{alg}\label{alg:C1}
Given $x \in \cap_{q\in\Z}G_q$ and $u \in G$, the following
algorithm computes $\mu_x(u)$.
\begin{algorithmic}\upshape
  \STATE Set $l=\len x$ and choose a permutation $q_0,q_1,\dots,q_l$ of the integers from $\inf x$ to $\sup x$.
  \STATE Set $u_0=u$.
  \FOR {$i=0$ to $l$ with step $+1$}
    \STATE Compute $u_i, \tilde\pi_{x,q_i}(u_i),
    \tilde\pi_{x,q_i}^2(u_i), \dots, \tilde\pi_{x,q_i}^N(u_i)$
    until repetition encountered.
    \STATE Set $u_{i+1}=\tilde\pi_{x,q_i}^N(u_i)$.
  \ENDFOR
  \STATE Set $v_{l+1}=u_{l+1}$.
  \FOR {$i=l$ to $0$ with step $-1$}
     \STATE Compute $v_{i+1}, \tilde\phi_{x,q_i}(v_{i+1}),
     \tilde\phi_{x,q_i}^2(v_{i+1}), \dots,
     \tilde\phi_{x,q_i}^N(v_{i+1})$ until repetition
     encountered and $u_i \prec \tilde\phi_{x,q_i}^N(v_{i+1})$.
     \STATE Set $v_i=\tilde\phi_{x,q_i}^N(v_{i+1})$.
  \ENDFOR
  \STATE \textbf{return} $v_0$
\end{algorithmic}
\end{alg}

\begin{proof}\label{alg:C2}
First, remark that $u_i \prec v_i$. By Proposition \ref{prop:B2},
there exists arbitrarily large $N$ such that $u_{i-1} \prec
\tilde\phi_{x,q_{i-1}}^N(v_i)$, so the algorithm stops in finite
steps. Moreover, we have $\Phi_xu_i \prec \Phi_xv_i$. Also notice
that, by Lemma \ref{lem:C2}, $\Phi_xv_0 \prec \Phi_xv_1 \prec
\cdots \prec \Phi_xv_{l+1} = \Phi_xu_{l+1} \prec \cdots \prec
\Phi_xu_1 \prec \Phi_xu_0$. Summarizing, we have $u_0 \prec v_0$
and $\Phi_xu_0 = \Phi_xv_0$. Therefore, by Lemma \ref{lem:C1}(4),
$\mu_x(v_0)=\mu_x(u_0)$.

Further, remark that $x^{v_i} \in G_{q_i}$ by Proposition
\ref{prop:B1}. By Lemma \ref{lem:B1}(1) and the latter claim of
Theorem \ref{thm:A1}, $x^{v_0} \in G_q$ for $\inf x \leq q \leq
\sup x$. Since $\inf x=\infs x$ and $\sups x = \sup x$, it follows
from the definition of $C^*(x)$ that $x^{v_0} \in C^*(x)$. Hence
$v_0=\mu_x(v_0)$.

Finally, we conclude that $v_0=\mu_x(v_0)=\mu_x(u_0)=\mu_x(u)$.
\end{proof}

Remark the inclusion $A(y) \subset \{ \mu_y(a) \mid a \in A \}$.
With above algorithm one may implement Algorithm \ref{alg:A3} by
computing the superset $\{ T(y^{\mu_y(a)}) \mid a \in A \}$
instead of $\{ T(y^u) \mid u \in A(y) \}$, both having a
cardinality not greater than the number of atoms of $G$.

Moreover, as in \cite{FG}, short-cuts can be used to increase the
efficiency. Actually, the set computed by the following algorithm
suffices for implementing Algorithm \ref{alg:A3}.

\begin{alg}
Given $x \in \cap_{q\in\Z}G_q$, the following algorithm computes a
set $\T$ satisfying $\{ T(x^u) \mid u \in A(x) \} \subset \T
\subset \{ T(x^{\mu_x(a)}) \mid a \in A \}$.
\begin{algorithmic}\upshape
  \STATE Set $Q=A$ and $\T=\emptyset$.
  \FOR {$a \in A$}
    \STATE Compute $\mu_x(a)$ by using Algorithm \ref{alg:C1}.
    \IF {meanwhile $a' \prec \tilde\pi_{x,q_i}^k(u_i)$ or $a' \prec \tilde\phi_{x,q_i}^k(v_i)$ for some $a' \in Q \setminus \{a\}$}
      \STATE Set $Q=Q\setminus\{a\}$.
    \ELSE
      \STATE Set $\T = \T \cup \{ T(x^{\mu_x(a)}) \}$.
    \ENDIF
  \ENDFOR
  \STATE \textbf{return} $\T$
\end{algorithmic}
\end{alg}

\begin{proof}
It is clear that $\T \subset \{ T(x^{\mu_x(a)}) \mid a \in A \}$.
Suppose an atom $a$ is excluded from $Q$ by another atom, say
$a_1$. Consider the sequence of atoms $a_1, a_2, \dots, a_k$ in
which $a_i$ is excluded from $Q$ by $a_{i+1}$ and $a_k$ survives
in $Q$ when the algorithm stops. By Lemma \ref{lem:C2}, $\Phi_xa_1
\prec \Phi_xa, \; \Phi_xa_2 \prec \Phi_xa_1, \dots, \Phi_xa_k
\prec \Phi_xa_{k-1}$ hence $\Phi_xa_k \prec \Phi_xa$. So, whenever
$\mu_x(a) \in A(x)$ we have $\Phi_xa_k = \Phi_xa$ and by Lemma
\ref{lem:C1}(1) $T(x^{\mu_x(a)}) = T(x^{\mu_x(a_k)}) \in \T$. The
inclusion $\{ T(x^u) \mid u \in A(x) \} \subset \T$ follows.
\end{proof}

It remains a proposition which is useful to compute the
pushforward and pullback of a simple element by means of simple
element calculus.

\begin{prop}
Given $x = \Delta^p x_1 \cdots x_l$ in normal form, $0 \leq k \leq
l$ and $u \in S \setminus \{\Delta\}$, define
$$\begin{array}{llll}
  u_0 = \tau^p(u), & u_i = \Delta \wedge x_i^{-1}u_{i-1}\Delta & \text{for} & i=1,\dots,k, \\
  u_{l+1} = u, & u_i = \Delta \wedge x_iu_{i+1} & \text{for} & i=l,\dots,k+1, \\
  v_k = u, & v_{i-1} = 1 \vee x_iv_i\Delta^{-1} & \text{for} & i=k,\dots,1, \\
  v_{k+1} = u, & v_{i+1} = 1 \vee x_i^{-1}v_i & \text{for} & i=k+1,\dots,l. \\
\end{array}
$$
Then $\phi_{x,p+k}(u) = u_k \wedge u_{k+1}$ and $\pi_{x,p+k}(u) =
\tau^{-p}(v_0) \vee v_{l+1}$.
\end{prop}

\begin{proof}
By induction one verifies the followings
$$\begin{array}{llll}
  u_i = \Delta \wedge x_i^{-1} \cdots x_1^{-1} \tau^p(u) \Delta^i & \text{for} & i=0,\dots,k, \\
  u_i = \Delta \wedge x_i \cdots x_l u & \text{for} & i=l+1,\dots,k+1, \\
  v_i = 1 \vee x_{i+1} \cdots x_k u \Delta^{i-k} & \text{for} & i=k,\dots,0, \\
  v_i = 1 \vee x_{i-1}^{-1} \cdots x_{k+1}^{-1} u & \text{for} & i=k+1,\dots,l+1. \\
\end{array}
$$
Let $x' = \Delta^px_1 \cdots x_k$ and $x''=x_{k+1} \cdots x_l$.
Note that $\phi_{x,p+k}(u) \in S$ and $\inf u = 0$. So
\begin{align*}
  & u_k \wedge u_{k+1}
  = \Big( \Delta \wedge x_k^{-1} \cdots x_1^{-1} \tau^p(u) \Delta^k \Big)
  \wedge \Big( \Delta \wedge x_{k+1} \cdots x_l u \Big)
  \\ & \quad \quad \quad \quad
  = \Delta \wedge x'^{-1}u\Delta^{p+k} \wedge x''u
  = \phi_{x,p+k}(u) \\
  & \tau^{-p}(v_0) \vee v_{l+1}
  = \tau^{-p}\Big( 1 \vee x_1 \cdots x_k u \Delta^{-k} \Big)
  \vee \Big( 1 \vee x_l^{-1} \cdots x_{k+1}^{-1} u \Big)
  \\ & \quad \quad \quad \quad
  = 1 \vee x'u\Delta^{-p-k} \vee x''^{-1}u
  = \pi_{x,p+k}(u)
\end{align*}
\end{proof}

\section{Revisiting braid groups}\label{sec:compute}

In this section braid groups are supposed to be endowed with the
classical Garside structure.

First, we argue that the algorithm by computing the ultra summit
set practically fails for solving conjugacy problem of reducible
braids. Let us consider a simple example. For any braid $\beta \in
B_{n-1}$ with $\inf\beta > 0$, appending one additional trivial
strand yields a reducible braid $\beta' \in B_n$. Note that the
cycling operation on $\beta'$ is essentially trivial (merely the
conjugation by the Garside element of $B_{n-1}$ on the subbraid
$\beta$), so the cycling-recurrence condition is always satisfied.
In the sequel, whenever $\beta$ lies in its super summit set, so
does $\beta'$ in its ultra summit set. Therefore, the ultra summit
set of $\beta'$ is at least as large as the super summit set of
$\beta$. As argued in \cite{Geb}, computation of super summit set
has been practically inaccessible for those braids with moderate
number of strands and word length, thus so is the computation of
ultra summit sets for such reducible braids.

However, one notices that the cycling operation $\cyc$ on the
components of a reducible braid can be achieved by applying
general cycling operations on the total braid. In above example,
to apply the cycling operation $\cyc$ on the subbraid $\beta$ it
suffices to apply $\cyc_{\inf\beta+1}$ on the total braid
$\beta'$. Then, with the cycling-recurrence condition posed on the
subbraid $\beta$, the cardinality $|C^*(\beta')|$ is comparable to
$|C^u(\beta)|$, contrasting with the fact that $|C^u(\beta')|$ is
not smaller than $|C^s(\beta)|$. Therefore, the set
$|C^*(\beta')|$ can still be effectively computed and the
conjugacy problem can be practically solved.

\smallskip

In the remainder of this section, we present some experimental
data to compare the performance of the ultra summit set $C^u$ with
the new summit set $C^*$ on solving conjugacy problem in braid
groups. In all tables, each entry involves a computation of 5,000
random braids.

\begin{table}
\caption{Experimental data for Test 1. Average/maximal sizes of
$C^u, C^*$ and average/maximal times $T^u, T^*$ spent on computing
them. 1K=1,000. Times are given in ms, unless stated otherwise.}
\label{tab:tab1}
\begin{center}
\begin{tabular}{c||c|c|c|c|c|c}
  $n$   & \multicolumn{6}{c}{5} \\
  \hline
  $l$       & 3         & 5         & 10        & 20        & 30        & 40 \\
  \hline
  $|C^u|$   & 21.6/48   & 81.4/168  & 599/3000  & 2345/64K  & 3760/239K & 4938/191K \\
  $|C^*|$   & 11.9/32   & 15.9/80   & 25.2/160  & 43.0/220  & 60.9/216  & 81.0/172 \\
  $T^u$     & 1/54      & 2/54      & 26/164    & 166/4560  & 364/23s   & 603/23s \\
  $T^*$     & 1/54      & 1/54      & 2/54      & 17/109    & 68/164    & 186/384 \\
  \multicolumn{7}{c}{} \\
  $n$   & \multicolumn{6}{c}{7} \\
  \hline
  $l$       & 3         & 5         & 10        & 20        & 30        & 40 \\
  \hline
  $|C^u|$   & 245/1824  & 7228/119K & ---       & ---       & ---       & --- \\
  $|C^*|$   & 31.3/288  & 27.2/612  & 33.7/352  & 59.0/176  & 87.5/164  & 117/152 \\
  $T^u$     & 11/109    & 443/7472  & ---       & ---       & ---       & --- \\
  $T^*$     & 1/54      & 1/54      & 5/109     & 43/164    & 184/439   & 529/769 \\
  \multicolumn{7}{c}{} \\
  $n$   & \multicolumn{6}{c}{9} \\
  \hline
  $l$       & 3         & 5         & 10        & 20        & 30        & 40 \\
  \hline
  $|C^u|$   & 3676/188K & ---       & ---       & ---       & ---       & --- \\
  $|C^*|$   & 41.3/1320 & 29.4/528  & 37.3/168  & 68.0/148  & 95.4/180  & 135/160 \\
  $T^u$     & 314/17s   & ---       & ---       & ---       & ---       & --- \\
  $T^*$     & 1/54      & 2/54      & 10/54     & 79/274    & 317/769   & 883/1318 \\
\end{tabular}
\end{center}
\end{table}

\begin{table}
\caption{Experimental data for Test 2. Average/maximal sizes of
$C^u, C^*$ and average/maximal times $T^u, T^*$ spent on computing
them. 1K=1,000. Times are given in ms, unless stated otherwise.}
\label{tab:tab2}
\begin{center}
\begin{tabular}{c||c|c|c|c|c|c}
  $n$   & \multicolumn{6}{c}{9} \\
  \hline
  $l$       & 2         & 3         & 5         & 10        & 15        & 20 \\
  \hline
  $|C^u|$   & 192/2740  & 267/6640  & 2681/197K & ---       & ---       & --- \\
  $|C^*|$   & 192/2740  & 66.4/2740 & 123/2880  & 416/14K   & 1070/38K  & 2770/125K \\
  $T^u$     & 10/164    & 17/604    & 225/24s   & ---       & ---       & --- \\
  $T^*$     & 4/109     & 1/54      & 3/109     & 29/879    & 160/11s   & 913/231s \\
  \multicolumn{7}{c}{} \\
  $n$   & \multicolumn{6}{c}{12} \\
  \hline
  $l$       & 2         & 3         & 5         & 10        & 15        & 20 \\
  \hline
  $|C^u|$   & 6064/355K & ---       & ---       & ---       & ---       & --- \\
  $|C^*|$   & 6064/355K & 445/77K   & 614/20K   & 3124/161K & 4121/80K  & 18K/440K \\
  $T^u$     & 641/38s   & ---       & ---       & ---       & ---       & --- \\
  $T^*$     & 252/15s   & 18/3021   & 29/824    & 406/26s   & 1346/98s  & 14s/979s \\
\end{tabular}
\end{center}
\end{table}

\begin{table}
\caption{Experimental data for Test 3. Average/maximal sizes of
$C^u, C^*$ and average/maximal times $T^u, T^*$ spent on computing
them. Times are given in ms, unless stated otherwise.}
\label{tab:tab3}
\begin{center}
\begin{tabular}{c||c|c|c|c|c|c}
  $n$   & \multicolumn{6}{c}{20} \\
  \hline
  $l$       & 5         & 10        & 20        & 30        & 40        & 50 \\
  \hline
  $|C^u|$   & 12.1/110  & 20.2/80   & 40.0/40   & 60.0/60   & 80.0/80   & 100.0/100 \\
  $|C^*|$   & 12.1/110  & 20.2/80   & 40.0/40   & 60.0/60   & 80.0/80   & 100.0/100 \\
  $T^u$     & 5/54      & 11/164    & 31/109    & 66/274    & 113/439   & 176/604 \\
  $T^*$     & 5/109     & 24/164    & 172/659   & 542/1538  & 1593/5109 & 3447/9780 \\
  \multicolumn{7}{c}{} \\
  $n$   & \multicolumn{6}{c}{50} \\
  \hline
  $l$       & 5         & 10        & 20        & 30        & 40        & 50 \\
  \hline
  $|C^u|$   & 10.0/20   & 20.0/20   & 40.0/40   & 60.0/60   & 80.0/80   & 100.0/100 \\
  $|C^*|$   & 10.0/20   & 20.0/20   & 40.0/40   & 60.0/60   & 80.0/80   & 100.0/100 \\
  $T^u$     & 24/109    & 38/109    & 106/329   & 208/604   & 351/1428  & 526/2307 \\
  $T^*$     & 19/54     & 78/219    & 518/1703  & 1544/4505 & 4461/12s  & 9609/23s \\
\end{tabular}
\end{center}
\end{table}

\begin{test}
This test compares the performance of $C^u$ and $C^*$ on the
reducible braids described in above example. For several values of
$n$ and $l$, we choose at random positive braids $\beta \in
B_{n-1}$ with $\sups\beta = l$. Then, for each of them we append
one additional trivial strand to make it into a reducible braid in
$B_n$ and compute the summit sets $C^u$ and $C^*$ of the braid
resulted. See Table \ref{tab:tab1}.

Random braids are generated as follows. Choose independent random
simple elements $x_1,x_2,\dots$ until $\sup (x_1 \cdots x_k) = l$.
Set $\beta=x_1 \cdots x_k$. Repeat this process until $\beta$
satisfies $\sups\beta = l$.
\end{test}

\begin{test}
This test compares the performance on a type of nested braids. For
several values of $n$ and $l$ with $n$ a multiple of three, we
choose at random positive braids $\beta \in B_3$ with $\sups\beta
= l$ in the same way as previous test. Then, for each $\beta$ we
choose independent random simple elements
$x_{i1},x_{i2},\dots,x_{il} \in B_{n/3}$ for $i=1,2,3$ and replace
each strand of $\beta$ by the braid $x_{i1} x_{i2} \cdots x_{il}$
to produce a nested braid of $n$ strands, then compute its summit
sets. See Table \ref{tab:tab2}.
\end{test}

\begin{test}
This test compares the performance on generic braids. For several
values of $n$ and $l$, we choose at random positive braids $\beta
\in B_n$ with $\lens\beta = l$ and compute the summit sets $C^u$
and $C^*$. See Table \ref{tab:tab3}.

Random braids are generated in the same way as \cite{Geb}. Choose
at random an integer $p \in \{0,1\}$ and choose independent random
simple elements $x_1,x_2,\dots$ until $\len (x_1 \cdots x_k) = l$.
Set $\beta=\Delta^px_1 \cdots x_k$. Repeat this process until
$\beta$ satisfies $\lens\beta = l$.
\end{test}

From above example and experimental data, we conclude that, with a
slight loss of efficiency for generic braids (in the worst case,
running time is prolonged approximately $\lens \beta$ times for a
braid $\beta$), a considerable improvement is achieved on solving
conjugacy problem of reducible braids by computing the new summit
set $C^*$ instead of the ultra summit set $C^u$.

\begin{rem}
Although the new summit set $C^*$ is very likely bounded above by
a polynomial function of word length for fixed number of strands
(see also \cite{BKL1,FG} for the conjectures on the bound of super
summit set), it is exponential in the number of strands. For
example, fix a braid $\beta \in B_3$ with $|C^*(\beta)|>1$, then
the new summit set $C^*$ of the reducible braid of $3n$ strands
yielded by juxtaposing $n$ copies of $\beta$ has a cardinality not
smaller than $|C^*(\beta)|^n$ which is obviously exponential in
$n$.

So, along the lines of solving conjugacy problem in braid groups
by computing some type of summit set, a polynomial algorithm both
in number of strands and word length will inevitably involve a
reduction process of reducible braids. Reader is referred to
\cite{BGG} for recent work on this direction. See also \cite{BNG,
Lee} for efforts on the relation between the reduction systems and
the super/ultra summit sets of reducible braids.
\end{rem}

\section{Other Applications}\label{sec:application}

In this section we present several more applications of the
machinery developed in the previous sections. More precisely, we
give simple proves, but in strong forms, to several results in
\cite{BKL2, LL1, BGG}.

\subsection{Summit infimum and supremum}\label{subsec:length}

The following main theorem of \cite{BKL2} gives rise to a bound
for computing summit infimum and supremum of a braid conjugacy
class by applying cycling and decycling operations. Indeed the
argument involved there is applicable to all Garside groups.
Recall that $\|\Delta\|$ denotes the norm of the Garside element,
which only depends on the Garside structure of a Garside group.

\begin{thm}[{\cite[Theorem 1]{BKL2}}]\label{thm:D0'}
If $\inf x < \infs x$ then $\inf\cyc^{\|\Delta\|-1}(x) > \inf x$.
Similarly, if $\sup x > \sups x$ then $\sup\dec^{\|\Delta\|-1}(x)
< \sup x$.
\end{thm}

Below we give a simple proof to a stronger version of above
theorem. (Note that the case $q=\inf x+1$ or $q=\sup x-1$ recovers
above theorem.)

\begin{thm}\label{thm:D0}
If $\inf x < q \leq \infs x$ then $\inf\cyc_q^{\|\Delta\|-1}(x)
\geq q$. Similarly, if $\sups x \leq q < \sup x$ then
$\sup\cyc_q^{\|\Delta\|-1}(x) \leq q$.
\end{thm}

\begin{proof}
We prove the first claim of the theorem and the latter claim can
be proved in the same way. Suppose $\inf x < q \leq \infs x$. By
Theorem \ref{thm:A0}, there exists $N>0$ such that $\inf
\cyc_q^N(x) \geq q$. Let $N$ be the minimum in possible. We have
to show that $N<\|\Delta\|$.

Set $y = (\tau^{-q}\cyc_q)^N(x) = \tau^{-Nq}\cyc_q^N(x)$ and
$\psi=\tau^{-q}\phi_{y,q}$. Then $\inf y \geq q$. Choose $u$ such
that $x = y^u$ and set $u' = \psi^N(u)^{-1}u$. Applying Lemma
\ref{lem:B3}(1) on the following diagram we obtain $x = y^{u'}$
and $\psi^n(u')=\psi^N(u)^{-1}\psi^n(u)$.

$$\begin{CD}
    x @>>> \tau^{-q}\cyc_q(x) @>>> \cdots @>>> (\tau^{-q}\cyc_q)^N(x) \\
    @AA{u}A @AA{\psi(u)}A @. @AA{\psi^N(u)}A \\
    (\tau^{-q}\cyc_q)^N(x) @= (\tau^{-q}\cyc_q)^N(x) @= \cdots @= (\tau^{-q}\cyc_q)^N(x) \\
    @AA{\psi^N(u)^{-1}}A @AA{\psi^N(u)^{-1}}A @. @AA{\psi^N(u)^{-1}}A \\
    (\tau^{-q}\cyc_q)^N(x) @= (\tau^{-q}\cyc_q)^N(x) @= \cdots @= (\tau^{-q}\cyc_q)^N(x) \\
    \\
  \end{CD}
$$
Moreover, by Lemma \ref{lem:B3}(2) and the minimality of $N$ we
have
$$1 = \psi^{N}(u') \precneqq \psi^{N-1}(u') \precneqq \cdots \precneqq \psi(u') \precneqq u'.$$

Note that by Lemma \ref{lem:B1}(2),(5) we have
$\psi^n(u'\wedge\Delta) = \psi^n(u')\wedge\Delta$. We claim that
$\inf y^{\psi^n(u'\wedge\Delta)} < q$ for $0 \leq n < N$.
Otherwise, by Lemma \ref{lem:B3}(2)
$$1 = \psi^N(u'\wedge\Delta) = \psi^{N-1}(u'\wedge\Delta) = \cdots = \psi^n(u'\wedge\Delta).$$
Hence we derive $\psi^n(u')\wedge\Delta = 1$, i.e. $\psi^n(u') =
1$, a contradiction.

Finally, applying Lemma \ref{lem:B3}(2) again yields
$$1 = \psi^N(u'\wedge\Delta) \precneqq \psi^{N-1}(u'\wedge\Delta)
  \precneqq \cdots \precneqq \psi(u'\wedge\Delta)
  \precneqq u'\wedge\Delta \precneqq \Delta
$$
which implies $N<\|\Delta\|$.
\end{proof}

\subsection{Stable summit set}\label{subsec:stable}

Very recently, the behavior of the cycling and decycling
operations on the powers of an element attracted much attention in
the study of Garside groups. See, for example, \cite{BGG, LL1,
LL2}. From this point of view, one perhaps is willing to introduce
the cycling operation of order $(p,q)$
$$\cyc_{p,q}(x) = x^{x^p \wedge \Delta^q}$$
and the $\cyc_{p,q}$-recurrence set
$$G_{p,q} = \{ x \in G \mid \text{$\cyc_{p,q}^N(x) = x$ for some $N>0$} \}.$$
Notice that $\cyc_q(x^p) = (\cyc_{p,q}(x))^p$, so applying a
$\cyc_q$ operation on $x^p$ is equivalent to applying a
$\cyc_{p,q}$ operation on $x$. In particular, $x^p \in G_q$ if and
only if $x \in G_{p,q}$.

Most arguments and results concerning the single order cycling
operations in this article can be generalized to the double order
version straightforwardly. The pushforward and pullback along the
cycling operation $x \to \cyc_{p,q}(x)$ are defined as
\begin{align*}
  & \phi_{x,p,q}(u) = x''u \wedge x'^{-1}\Delta^q\tau^q(u), \\
  & \pi_{x,p,q}(u) = \Delta^{\inf u} \vee x''^{-1}u \vee x'\Delta^{-q}\tau^{-q}(u),
\end{align*}
respectively, where $x'=x^p\wedge\Delta^q$ and $x''=x'^{-1}x^p$.

With a suitable modification, the algorithms for computing
$C^*(x)$ can be used to compute the set
$$C^{[m,n],*}(x) = C(x) \cap \bigcap_{m \leq p \leq n, \; q \in \Z} G_{p,q}$$
which is hence nonempty.

However, it should be pointed out that the knowledge of present
article does not lead to an algorithm to compute the refined
summit set subject to all $\cyc_{p,q}$-recurrence conditions
$$C^{*,*}(x) = C(x) \cap \bigcap_{p,q \in \Z} G_{p,q},$$
because we do not know how to bound the order $p$. Nevertheless we
have the following theorem.

\begin{thm}\label{thm:D1}
The set $C^{*,*}(x)$ is nonempty.
\end{thm}

\begin{proof}
Note that $C^{*,*}(x)$ is the intersection of the descending
sequence of finite, nonempty sets
$$C^{[-1,1],*}(x) \supset C^{[-2,2],*}(x) \supset C^{[-3,3],*}(x) \supset \cdots.$$
Since a descending sequence of finite, nonempty sets always has
nonempty intersection (we leave it to the reader as an easy
excise), the theorem follows.
\end{proof}

In the sequel, as supersets of $C^{*,*}(x)$ the {\em stable super
summit set} \cite{LL1}
$$C(x) \cap \bigcap_{p\geq1, \; q \in \{\infs x^p,\sups x^p\}} G_{p,q}$$
and the {\em stable ultra summit set} \cite{BGG}
$$C(x) \cap \bigcap_{p\geq1, \; q \in \{\infs x^p,\infs x^p+1,\sups x^p\}} G_{p,q}$$
are both nonempty.

\subsection{Rigid elements}\label{subsec:rigid}

Rigid elements became of interest in the conjugacy problem in
Garside groups because, on the one hand, these elements have many
nice properties and, on the other hand, these elements are generic
enough, for example, it is shown in \cite{BGG} that for every
pseudo-Anosov braid $x$ some power of $x$ is rigid up to
conjugacy.

In \cite{BGG}, an element $x = \Delta^p x_1 \cdots x_l$ in normal
form is said to be {\em rigid} if $l>0$ and
$$\Delta^p x_1 \cdots x_l \tau^{-p}(x_1)$$
is also in normal form. Actually the second condition is
equivalent to say
$$x^2 = \Delta^{2p} \tau^p(x_1) \cdots \tau^p(x_l) x_1 \cdots x_l$$
is in normal form. So the condition in the definitions can be
stated more intrinsically as $\len x>0$ and
$$x^2 \wedge \Delta^{\inf x+\sup x} = x\Delta^{\inf x}.$$

As mentioned above, rigid elements have many nice properties. For
example, a nontrivial power of a rigid element is also rigid.
Another example is the behavior of the cycling operations on them
is very simple. Indeed, the action of each $\cyc_{p,q}$ on a rigid
element $x = \Delta^p x_1 \cdots x_l$ in normal form is merely a
cyclic permutation together with some possible $\tau$ actions on
the $x_i$'s. It follows that the rigid elements of a Garside group
$G$ are contained in $\cap_{p,q\in\Z}G_{p,q}$

The following theorem is one of the main result in \cite{BGG}.

\begin{thm}[{\cite[Theorem 3.22 and Theorem 3.34]{BGG}}]\label{thm:D2'}
(1) If $\lens x>1$ and $x$ is rigid then the ultra summit set
$C^u(x)$ is precisely the rigid conjugates of $x$.

(2) If $\lens x>1$ and a power $x^N$ is conjugate to a rigid
element, then $N$ can be chosen so that $0 < N < \|\Delta\|^2$.
\end{thm}

Below we give an alternative proof to above theorem, with the
boring condition $\lens x>1$ dropped.

\begin{thm}\label{thm:D2}
(1) If $x$ is rigid then the set $C^{*,*}(x)$ is precisely the
rigid conjugates of $x$.

(2) If a power $x^N$ is conjugate to a rigid element, then $N$ can
be chosen so that $0 < N < \|\Delta\|^2$.
\end{thm}

The proof depends on several lemmas.

\begin{lem}\label{lem:D21}
Suppose $x$ is rigid and $p=\inf x$, $q=\inf x+\sup x$. Then for
$x^u \in C^s(x)$ we have $u \prec \tau^{-p}\phi_{x,2,q}(u)$ with
equality holds if and only if $x^u$ is rigid.
\end{lem}

\begin{proof}
By the definition of rigidity, $x^2\wedge\Delta^q = x\Delta^p$. A
straightforward calculation shows that $\tau^{-p}\phi_{x,2,q}(u) =
x^{-1} u ((x^u)^2\wedge\Delta^q) \Delta^{-p}$ (see the diagram
below).
$$\begin{CD}
    x^u @>{(x^u)^2\wedge\Delta^q}>> \cyc_{2,q}(x^u) @>{\Delta^{-p}}>> \tau^{-p}\cyc_{2,q}(x^u) \\
    @A{u}AA @AA{\phi_{x,2,q}(u)}A @AA{\tau^{-p}\phi_{x,2,q}(u)}A \\
    x @>{x\Delta^p}>> \cyc_{2,q}(x) @>{\Delta^{-p}}>> x \\
    \\
  \end{CD}
$$
Since both $x,x^u \in C^s(x)$, we have $\len x^u = \len x > 0$ and
$\inf x^u=p$, $\sup x^u=q-p$. Therefore,
$$u \prec u (x^u\Delta^{-p} \wedge (x^u)^{-1}\Delta^{q-p})
  = x^{-1} u ((x^u)^2\wedge\Delta^q) \Delta^{-p}
  = \tau^{-q}\phi_{x,2,q}(u)
$$
with equality holds if and only if $x^u\Delta^{-p} \wedge
(x^u)^{-1}\Delta^{q-p} = 1$, i.e. $(x^u)^2\wedge\Delta^q =
x^u\Delta^p$, that is, $x^u$ is rigid.
\end{proof}

\begin{lem}\label{lem:D22}
If $x^m$ $(m>0)$ is rigid and $\inf x^m=m\inf x$, $\sup x^m=m\sup
x$ then $x$ is also rigid.
\end{lem}

\begin{proof}
Clearly $\len x = \len x^m/m > 0$. Let $p=\inf x$, $q=\inf x+\sup
x$. Then $\inf x^m = mp$, $\sup x^m = mq-mp$ and the rigidity of
$x^m$ says $x^{2m}\wedge\Delta^{mq} = x^m\Delta^{mp}$. Therefore,
\begin{align*}
  x\Delta^p
  & \prec x^2 \wedge \Delta^q
  \prec x^{(m+1)}\Delta^{-(m-1)p} \wedge x^{-(m-1)}\Delta^{-(m-1)p+mq} \\
  & = x^{-(m-1)} (x^{2m}\wedge\Delta^{mq}) \Delta^{-(m-1)p}
  = x^{-(m-1)} (x^m\Delta^{mp}) \Delta^{-(m-1)p} = x\Delta^p,
\end{align*}
which implies $x^2\wedge\Delta^q = x\Delta^p$. So, $x$ is rigid.
\end{proof}

\begin{lem}\label{lem:D23}
For every $x$ there exist integers $0 < N_1,N_2 \leq \|\Delta\|$
such that $\infs x^{nN_1} = n\infs x^{N_1}$ and $\sups x^{nN_2} =
n\sups x^{N_2}$ for all $n>0$.
\end{lem}

\begin{proof}
Let $N_1,N_2$ be the denominators of the rational numbers
$$\max \{\infs x^n/n \mid n=1,2,\dots,\|\Delta\| \}$$
and
$$\min \{\sups x^n/n \mid n=1,2,\dots,\|\Delta\| \}$$
respectively. Clearly, $0 < N_1,N_2 \leq \|\Delta\|$. Then
\cite[Lamma 2.4(ii), Theorem 3.2 and Theorem 5.1(i)]{LL2} says
$N_1,N_2$ are exactly what we want.
\end{proof}

\begin{proof}[Proof of Theorem \ref{thm:D2}]
(1) Let $p=\inf x$, $q=\inf x+\sup x$. Since the rigid elements of
$G$ are contained in $\cap_{p,q\in\Z}G_{p,q}$, the rigid
conjugates of $x$ belong to $C^{*,*}(x)$. Moreover, we have the
obvious inclusion $C^{*,*}(x) \subset C^s(x) \cap G_{2,q}$. So, it
suffices to show that all elements of $C^s(x) \cap G_{2,q}$ are
rigid.

Suppose $x^u \in C^s(x) \cap G_{2,q}$. By Lemma \ref{lem:D21},
$$u \prec \tau^{-p}\phi_{x,2,q}(u) \prec \cdots \prec (\tau^{-p}\phi_{x,2,q})^N(u).$$
Since both $x,x^u \in G_{2,q}$, by the same argument in the proof
of Proposition \ref{prop:B1} we have $(\tau^{-p}\phi_{x,2,q})^N(u)
= u$ for some $N>0$. Therefore, $u = \tau^{-p}\phi_{x,2,q}(u)$,
and Lemma \ref{lem:D21} says that $x^u$ is rigid.

(2) Suppose a power of $x$, say $x^m$ with $m>0$, is conjugate to
a rigid element. By Theorem \ref{thm:D1} we can choose $y \in
C^{*,*}(x)$. Then $y^n \in C^{*,*}(x^n)$ hence $\inf y^n = \infs
x^n$ and $\sup y^n = \sups x^n$ for all $n \in \Z$.

Let $N_1,N_2$ be as in Lemma \ref{lem:D23} and let $N$ be the
least common multiple of $N_1,N_2$. Note that $N<\|\Delta\|^2$;
otherwise, we must have $\|\Delta\|=1$, which implies $G$ is the
infinite cyclic group generated by $\Delta$ hence has no rigid
element.

Since $x^{mN}$ is also conjugate to a rigid element, $y^{mN} \in
C^{*,*}(x^{mN})$ is rigid by the conclusion of (1). Moreover, from
the choice of $N$ we have $\inf y^{mN} = m\inf y^N$ and $\sup
y^{mN} = m\sup y^N$. By Lemma \ref{lem:D22} $y^N$ is rigid. This
completes the proof of the theorem.
\end{proof}

From the proof we have the follow algorithm for deciding whether
some power of $x$ is conjugate to a rigid element and computing
the possible rigid element. First, compute the integers $0 <
N_1,N_2 \leq \|\Delta\|$ as in Lemma \ref{lem:D23} (which indeed
can be done in polynomial time in the case of braid groups). Then,
compute an element $\tilde{x}$ of $C^s(x^N) \cap G_{2,\infs
x^N+\sups x^N}$ where $N$ is the least common multiple of
$N_1,N_2$. Then some power of $x$ is conjugate to a rigid element
if and only if $\tilde{x}$ is rigid.

Moreover, the proof says for every rigid element $x$ we have
$$C^{*,*}(x) = C^s(x) \cap G_{2,\inf x+\sup x}$$
which is precisely the rigid conjugates of $x$.


\begin{thebibliography}{99}
\bibitem{BGG} J. S. Birman, V. Gebhardt, J. Gonz\'alez-Meneses,
    Conjugacy in Garside Groups I: Cyclings, Powers, and Rigidity,
    arXiv:math.GT/0605230.
\bibitem{BKL1} J. S. Birman, K. H. Ko, S. J. Lee,
    A new approach to the word and conjugacy problems in the braid groups,
    Adv. Math. 139(2) (1998), 322--353.
\bibitem{BKL2} J. S. Birman, K. H. Ko, S. J. Lee,
    The infimum, supremum, and geodesic length of a braid conjugacy class,
    Adv. Math. 164(1) (2001), 41--56.
\bibitem{BNG} D. Bernardete, Z. Nitecki, M. Guti\'errez,
    Braids and the Nielsen-Thurston classification,
    J. Knot Theory Ramif. 4(4) (1995), 549--618.
\bibitem{EM} E. A. El-Rifai, H. R. Morton,
    Algorithms for positive braids,
    Quart. J. Math. Oxford 45 (1994), 479--497.
\bibitem{Deh} P. Dehornoy,
    Groupes de Garside,
    Ann. Sci. \'Ecole Norm. Sup. (4) 35(2) (2002), 267--306.
\bibitem{DP} P. Dehornoy, L. Paris,
    Gaussian groups and Garside groups, two generalizations of Artin groups,
    Proc. London Math. Soc. 79(3) (1999), 569--604.
\bibitem{FG} N. Franco, J. Gonz\'alez-Meneses,
    Conjugacy problem for braid groups and Garside groups,
    J. Algebra 266(1) (2003), 112--132.
\bibitem{Gar} F. A. Garside,
    The braid group and other groups,
    Quart. J. Math. Oxford Ser. 20(2) (1969), 235--254.
\bibitem{Geb} V. Gebhardt,
    A new approach to the conjugacy problem in Garside groups,
    J. Algebra 292(1) (2005), 282--302.
\bibitem{Lee} S. J. Lee,
    Garside theory on reducible braids,
    arXiv:math.GT/0506188.
\bibitem{LL1} E. K. Lee, S. J. Lee,
    Stable super summit sets in Garside groups,
    arXiv:math.GT/0602582.
\bibitem{LL2} E. K. Lee, S. J. Lee,
    Some power of an element in a Garside group is conjugate to a periodically geodesic element,
    arXiv:math.GN/0604144.
\bibitem{Pic} M. Picartin,
    The conjugacy problem in small Gaussian groups,
    Comm. Algebra 29(3) (2001), 1021--1039.
\end{thebibliography}
\end{document}